\documentclass[12pt]{article}
\usepackage{graphicx,amssymb, amsmath}
\title{The six-dimensional Delaunay polytopes}

\newtheorem{proposition}{Proposition}
\newtheorem{theorem}{Theorem}
\newtheorem{corollary}{Corollary}

\newtheorem{definition}{Definition}
\newtheorem{remark}{Remark}
\newcommand{\proof}{\noindent{\bf Proof.}\ \ }

\author{Mathieu Dutour\\
ENS, Paris and Hebrew University, Jerusalem,\footnote{Research financed by EC's IHRP Programme, within the Research Training Network ``Algebraic Combinatorics in Europe,'' grant HPRN-CT-2001-00272.}
}

\begin{document}
\newcommand{\R}{\ensuremath{\mathbb{R}}}
\newcommand{\N}{\ensuremath{\mathbb{N}}}
\newcommand{\Q}{\ensuremath{\mathbb{Q}}}
\newcommand{\C}{\ensuremath{\mathbb{C}}}
\newcommand{\Z}{\ensuremath{\mathbb{Z}}}
\newcommand{\T}{\ensuremath{\mathbb{T}}}

\maketitle

\begin{abstract}
\noindent Given a lattice $L$, a full dimensional polytope $P$ is called a
{\em Delaunay polytope} if the set of its vertices is $S\cap L$ with $S$ being
an {\em empty sphere} of the lattice. Extending our previous work 
\cite{DD-hyp} on the {\em hypermetric cone} $HYP_7$, we classify the
six-dimensional Delaunay polytopes according to their {\em combinatorial 
type}. The list of $6241$ combinatorial types is obtained by a study of 
the set of faces of the polyhedral cone $HYP_7$.
\end{abstract}

\section{Introduction}
A {\em distance vector} $(d_{ij})_{0\leq i< j\leq n}\in \R^N$ with $N={n+1\choose 2}$ is called an {\em $(n+1)$-hypermetric} if it satisfies the following {\em hypermetric inequalities}:
\begin{equation}\label{Definition-of-hypermetrics}
H(b)d=\sum_{0\leq i<j\leq n} b_ib_jd_{ij}\leq 0\mbox{~for~any~}b=(b_i)_{0\leq i\leq n}\in\Z^{n+1}\mbox{~with~}\sum_{i=0}^{n}b_i=1\,\,.
\end{equation}
The set of distance vectors satisfying (\ref{Definition-of-hypermetrics}) is called the {\em hypermetric cone} and denoted by $HYP_{n+1}$.

In fact, $HYP_{n+1}$ is a polyhedral cone (see \cite{DGL93} p.~199). Lovasz (see \cite{DL} p.~201-205) gave another proof of it and bound $\max |b_i|\leq n!2^n{2n \choose n}^{-1}$ for any vector $b=(b_i)_{0\leq i\leq n-1}$ defining a facet of $HYP_n$.

There is a one-to-one correspondence between {\em non-degenerate} elements of 
$HYP_{n+1}$ and semi-metrics on {\em affine basis} of $n$-dimensional 
Delaunay polytopes. So, the enumeration of combinatorial types of 
$n$-dimensional Delaunay polytopes is reduced to the enumeration of 
non-degenerate faces of $HYP_{n+1}$ under an equivalence relation 
called {\em geometrical equivalence} (see Remark 
\ref{The-Simple-Case-Of-Geometrical-Equivalence} and Section 
\ref{Enumeration-Techniques}). Moreover, the dimension of faces 
of $HYP_{n+1}$ allows to define the notion of {\em rank} of a Delaunay 
polytope (see \cite{DGL92}).

We have the inclusion $CUT_n\subset HYP_n$, where $CUT_n$ (see Section \ref{Section-DimensionSix} below and Chapter $4$ of \cite{DL}) is the cone generated by all cut semi-metrics on $n$ points.
One has $HYP_n=CUT_n$ for $n\leq 6$; so, the enumeration of combinatorial types of Delaunay polytopes of dimension less or equal to five correspond to the study of faces of $CUT_6$. This study was done by Fedorov (\cite{Fe85}), Erdahl and Ryshkov (\cite{ErRyI} and \cite{ErRyII}), and Kononenko (\cite{Ko3}) in dimension three, four and five, respectively.

In the case of dimensional six, the inclusion $CUT_7\subset HYP_7$ is strict,
but the description of facets (i.e. faces of rank $20$, corresponding
to {\em repartitioning polytopes}) and extreme rays (i.e. faces of
rank $1$, corresponding to {\em extreme Delaunay polytopes})
of the cone $HYP_7$ is still possible (see \cite{DD-hyp}).
In Section \ref{Section-DimensionSix} we present the number of
combinatorial types of six-dimensional Delaunay polytope for every
given rank. In Section \ref{Enumeration-Techniques} we present our
method for the computation of those combinatorial types using the
face lattice of $HYP_7$.


Voronoi \cite{Vo} defines a partition of the cone $PSD_n$ of positive 
semi-definite quadratic forms by {\em $L$-type domains}. Two forms
of the same $L$-type domain have affinely equivalent Voronoi polytopes.
Every vertex of a Voronoi polytope correspond to a center of a Delaunay
polytope. All Delaunay polytopes form a partition of the space $\R^n$,
which is dual to the partition by Voronoi polytopes. It is proved in
\cite{DGL93} that the hypermetric cone $HYP_n$ is the union of a finite
number of $L$-type domains.

While Voronoi theory of $L$-type domains describes the combinatorial
structure of lattices, the theory of hypermetrics (\cite{DGL92}
and \cite{DL}) describes combinatorial structure of one 
Delaunay polytope in a lattice. While the hypermetric cone $HYP_n$ has
the symmetry group $Sym(n)$, a $L$-type domain, in general, has
trivial group.

\section{Delaunay polytopes and hypermetrics}\label{General-result}
\noindent Here we present the main notions for hypermetrics, which are needed to our study; the presentation is slightly simplified by the systematic use of affine bases. For the complete theory, with proofs, see \cite{DGL92} and Chapters $13$--$16$ of \cite{DL}.

A family $v_0, \dots, v_n$ of $n+1$ vertices of $\R^n$ is called {\em independent} if the family $(v_i-v_0)_{1\leq i\leq n}$ has linear rank $n$.
Let $L\subset \R^n$ be a $n$-dimensional lattice and let $S=S(c,r)$ be a  sphere in $\R^n$ with center $c$ and radius $r$. Then $S$ is said to be an {\em empty sphere} in $L$ if the following two conditions hold:
\begin{itemize}
\item[(i)] $\Vert v-c\Vert\geq r$ for all $v\in L$ and
\item[(ii)] the set $S\cap L$ contain an independent set of size $n+1$.
\end{itemize}
The center $c$ of $S$ is called a {\em hole} in \cite{CS}. The polytope $D$, which is defined as the convex hull of the set $S\cap L$ is called a {\em Delaunay polytope}, or (in original terms of Voronoi, who introduced them in \cite{Vo}) {\em L-polytope}.

\begin{definition}
Let $D$ be a $n$-dimensional Delaunay polytope with vertex-set $V$.

(i) A family  $v_0$, \dots, $v_n$ of vertices of $D$ is called an {\em affine basis} if for all $v\in V$ there exist an unique family $(b_i)_{0\leq i\leq n}\in \Z^{n+1}$, such that
\begin{equation*}
\sum_{i=0}^{n}b_i v_i=v\mbox{~and~}\sum_{i=0}^{n}b_i=1\;.
\end{equation*}

(ii) The Delaunay polytope $D$ is called {\em basic} if it has at least one affine basis. The vertices of an affine basis are called {\em basic vertices}.

\end{definition}
All known Delaunay polytopes are basic. We prove in 
Theorem \ref{ThanksSlava} that all six-dimensional Delaunay polytopes 
are basic too. We will always assume that Delaunay polytopes are basic.

For every family ${\cal A}=\{ v_0, \dots, v_m\}$ of vertices of a Delaunay polytope $P$ one can define a distance function $d_{\cal A}$ by $d_{\cal A}(i,j)=\Vert v_i-v_j\Vert^2$. The function $d_{\cal A}$ turns out to be an hypermetric by the following formula (see \cite{As} and \cite{DL} p.~195) :
\begin{equation*}
\sum_{0\leq i,j\leq m} b_ib_jd_{\cal A}(i,j)=2(r^2-\Vert \sum_{i=0}^m b_iv_i - c\Vert^2)\leq 0\,.
\end{equation*}
On the other hand, Assouad has shown in \cite{As} that {\em every} $d\in HYP_{n+1}$ can be expressed as $d_{\cal A}$ with ${\cal A}$ being a family of vertices of a Delaunay polytope $D$ of dimension {\em less or equal} to $n$.

A ray $d\in HYP_{n+1}$ is called {\it non-degenerate} if $d=d_{\cal B}$ with ${\cal B}$ an affine basis of a $n$-dimensional Delaunay polytope $D$.
For a given ray $d\in HYP_{n+1}$ the {\it annulator} is defined by

\begin{equation*}
Ann(d)=\{b\in\Z^{n+1}\mbox{~:~}\sum_{i=0}^n b_i=1\mbox{~and~}H(b)d=0\}.
\end{equation*}
We call {\em basic vectors} the vectors $e_i=(0,\dots, 0, 1, 0, \dots, 0)\in\Z^{n+1}$ with $0\leq i\leq n$. Using Proposition below, those vectors $e_i$ are identified with basic vertices $v_i$. One has $H(b)=0$ if and only if $b=e_i$ for some $0\leq i\leq n$.

\begin{proposition}\label{BasicReductionResult}
Let $D$ be an $n$-dimensional Delaunay polytope with vertex-set $V$; let $d=d_{\cal A}$ with ${\cal A}=\{v_0, \dots, v_{n}\}$ a family of $V$. Then there is equivalence between following properties:

(i) the family ${\cal A}$ is independent and

(ii) $det\,[d(0,i)+d(0,j)-d(i,j)]_{1\leq i,j\leq n}\not= 0$.

Also one has equivalence between properties:

(iii) ${\cal A}$ is an affine basis and

(iv) the following mapping is one-to-one
\begin{equation*}
b\in Ann(d)\mapsto \sum_{i=0}^n b_i v_i\in V \;.
\end{equation*}

\end{proposition}
\proof The family $(v_i)_{0\leq i\leq n}$ is independent if and only if the matrix $((v_i-v_0).(v_j-v_0))_{1\leq i,j\leq n}$ is positive definite. The formula $2(v_i-v_0).(v_j-v_0)=\Vert v_i-v_0\Vert^2+\Vert v_j-v_0\Vert^2-\Vert v_i-v_j\Vert^2$, gives the first equivalence. The second equivalence is obvious.

By above Proposition, the set $Ann(d)$ is finite if $d$ is non-degenerate, since the hypermetric vectors $b\in Ann(d)$ correspond to lattice points $v$, which belong to a sphere. We will consider only non-degenerate hypermetrics (in the hypermetric cone $HYP_{n+1}$) associated with affine bases of $n$-dimensional Delaunay polytopes.

Using above mapping, we identify vertices $v$ of $D$ with hypermetric vectors $b$. For example, for the affine basis ${\cal B}=\{v_0,\dots, v_3\}$ with
$v_0=(0,0,0)$, $v_1=(1,0,0)$, $v_2=(0,1,0)$ and $v_3=(0,0,1)$ of the $3$-cube,  one has $(1,1,0)=v_1+v_2-v_0$, $(0,1,1)=v_2+v_3-v_0$, $(1,0,1)=v_1+v_3-v_0$ and $(1,1,1)=v_1+v_2+v_3-2v_0$. 
Note that, while $(-2,1,1,1)\in Ann(d_{\cal B})$, the inequality $H(-2,1,1,1)d\geq 0$ does not define a facet of the cone $HYP_4$, since $H(-2,1,1,1)=H(-1,1,1,0)+H(-1,0,1,1)+H(-1,1,0,1)$.

The face $F(d)$, associated to an hypermetric $d\in HYP_{n+1}$, is the minimal face, containing the vector $d$. It can also be defined as:
\begin{equation*}
F(d)=\{e\in HYP_{n+1}\mbox{~:~}H(b)e=0\mbox{~for~all~}b\in Ann(d).\}
\end{equation*}

The {\em rank} of a $n$-dimensional Delaunay polytope $D$ is defined as 
dimension of $F(d_{\cal B})$, where ${\cal B}$ is an affine basis of $D$ 
(see \cite{DL} p.~217); its {\em corank} is, by definition,
${n+1\choose 2}-rank\,\,D$. The rank of a Delaunay polytope 
$D$ is equal to the topological dimension of the set of affine 
bijections $T$ of $\R^{n}$ (up to translations and orthogonal 
transformations), for which $T(D)$ is again a Delaunay polytope 
(see \cite{DL} p.~225).

If $D$ is an $n$-simplex, then its unique affine basis defines an hypermetric $d$, which has $Ann(d)=\{e_0, \dots, e_n\}$ and $F(d)=HYP_{n+1}$.
If $D$ is a $n$-dimensional Delaunay polytope with $n+2$ vertices $v_0,\dots, v_{n+1}$ (such polytope is called {\em repartitioning polytope}) and affine basis ${\cal B}=\{v_0,\dots, v_{n}\}$, then, by writing $v_{n+1}=\sum_{i=0}^n b_i v_i$, one obtains $F(d_{\cal B})=\{d\in HYP_{n+1}, H(b)d=0\}$. So, $F(d)$ is a facet of $HYP_{n+1}$. On the other hand, if $d\in HYP_{n+1}$ with $F(d)$ being a facet, then $Ann(d)=\{e_0, \dots, e_n, b\}$ and $d=d_{\cal B}$ with ${\cal B}$ being an affine basis of a repartitioning polytope.

Define two Delaunay polytopes to be {\em affinely equivalent}, if there is an affine bijective mapping transforming one into another; the equivalence classes by this relation are called {\em combinatorial type}. By Theorem 15.2.1 in \cite{DL} p.~222, two vectors $d$, $d'$ of $HYP_{n+1}$, such that $F(d)=F(d')$, correspond to affinely equivalent Delaunay polytopes. Moreover, if two Delaunay polytopes $D$, $D'$ are affinely equivalent, then an affine basis ${\cal B}$ of $D$ induces an affine basis ${\cal B'}$ of $D'$, such that $F(d_{\cal B})=F(d_{\cal B'})$.

A face of $HYP_{n+1}$ is called {\em non-degenerate}, if one of its interior vectors $d$ is non-degenerate (by Theorem 15.2.1 of \cite{DL}, this is equivalent to {\em all} interior vectors being non-degenerate).
One can consider only non-degenerate faces, since we already assumed, that all Delaunay polytopes are basic.
Given a face $F$ of $HYP_{n+1}$, we define 
$$Ann(F)=\{b\in\Z^{n+1}\mbox{~~:~~}\sum_{i}b_i=1\mbox{~and~}H(b)d=0\mbox{~for~all~}d\in F \}.$$
We will call $\{b^0, \dots, b^n\}\subset Ann(F)$ an {\em affine basis of the combinatorial type} if $det\,\{b^0, \dots, b^n\}=\pm 1$. Any affine basis of $Ann(F)$ corresponds to an affine basis of a Delaunay polytope $D$, if $D$ has this combinatorial type.

\begin{definition}\label{DefinitionGeometricalEquivalence}
Two faces $F$, $F'$ of $HYP_{n+1}$ are said to be {\em geometrically equivalent} if there exist an affine basis ${\cal B}=\{b^0, \dots, b^n\}$ of $Ann(F)$, such that the mapping
\begin{equation*}
\left\lbrace\begin{array}{rcl}
\phi_{F',F}:Ann(F')&\rightarrow&Ann(F)\\
b&\mapsto & b_0b^0+\dots b_nb^n
\end{array}\right.
\end{equation*}
is bijective.

\end{definition}

If $\pi$ is a permutation of $\{0,\dots, n\}$ and $F$ a face of $HYP_{n+1}$, then $\{e_{\pi(0)},\dots, e_{\pi(n)}\}$ is an affine basis of $Ann(F)$; so, two faces $F$ and $F'$, which are equivalent, by a symmetry of $HYP_{n+1}$, are also geometrically equivalent. The reverse is not true, in general.

\begin{remark}\label{The-Simple-Case-Of-Geometrical-Equivalence}
Let $D$ be a $n$-dimensional repartitioning polytope with an affine basis ${\cal B}=\{v_0, \dots, v_{n}\}$; we write $v_{n+1}=\sum_{i=0}^n b_i v_i$. The family ${\cal B'}=(v_0, \dots v_{i-1}, v_{n+1}, v_{i+1}, \dots, v_n)$ is an affine basis if and only if $|b_i|=1$. 

If $b_i=-1$, then $F(d_{\cal B'})=F(d_{\cal B})=\{d\in HYP_{n+1}, H(b)d=0\}$, while if $b_i=1$, then $F(d_{\cal B})\not= F(d_{\cal B'})=\{d\in HYP_{n+1}, H(b')d=0\}$ with $b'=(-b_0, \dots, -b_{i-1}, b_i, -b_{i-1}, \dots, -b_n)$.
\end{remark}

So, the study of combinatorial types of $n$-dimensional Delaunay polytopes is reduced to the study of non-degenerate faces of the hypermetric cone $HYP_{n+1}$ under the geometrical equivalence.

Consider, as an example, the two-dimensional case. We have $HYP_3=CUT_3$ and $CUT_3$ has three facets, which correspond to $H(-1,1,1)$ and its permutations. One can check that a triangle ${\cal T}=\{v_0,v_1,v_2\}$ satisfies $d_{\cal T}=(d_{01},d_{02},d_{12})\in HYP_3$ if and only if it is an obtuse triangle. 
The vector $d_{\cal T}$ is non-degenerate if and only if three vertices of ${\cal T}$ are not aligned. 
Moreover, $d_{\cal T}$ is incident to an hypermetric facet, say, $H(-1,1,1)$, if and only if the vertex $v_0$ has angle $\frac{\pi}{2}$, in which case the Delaunay polytope has four vertices: $v_0$, $v_1$, $v_2$ and $v_1+v_2-v_0$. 
So, there are two combinatorial types of Delaunay polytopes in dimension two: obtuse triangles and rectangles.

Take a face $F$ of the hypermetric cone $HYP_{n+1}$. To every
$b=(b_i)_{0\leq i\leq n}\in Ann(F)$, we associate a vertex
$(b_1, \dots, b_n)\in\Z^n$; let the set of such vertices being 
denoted by $V$. Every distance vector $d\in F$ correspond
to a Gram matrix $G$. The set $V$ is the vertex-set of a Delaunay
polytope, which is described by the scalar product defined by $G$.
So, we can encode, in our computations, the combinatorial type of
Delaunay polytope by the set $Ann(F)$.

\begin{proposition}
If $D$ is a Delaunay polytope and ${\cal B}$, ${\cal B'}$ are two affine bases of $D$, then the faces $F(d_{\cal B})$, $F(d_{\cal B'})$ are equivalent up to a linear mapping.
This linear mapping preserve the non-degeneracy.
\end{proposition}
\proof If $(v_0, \dots, v_{n})$ and $(v'_0, \dots, v'_{n})$ are two affine bases of $D$, then one can express $v'_i$ in terms of $(v_j)_{0\leq i\leq n}$ as follows
\begin{equation*}
v'_i=\sum_{j=0}^{n}\alpha_{ij}v_j\mbox{~~with~~}\sum_{j=0}^{n} \alpha_{ij}=1\;.
\end{equation*}
One can express $||v'_{i_1}-v'_{i_2}||^2$ in terms of $||v_{j_1}-v_{j_2}||^2$. This induces a linear mapping $\phi$ from $F(d_{\cal B})$ to $F(d_{\cal B'})$; 
expressing $v_j$ in terms of $(v'_i)_{0\leq i\leq n}$, one get the reverse mapping $\phi^{-1}$ and so, the linear equivalence.

If $d\in F(d_{\cal B})$, then $d$ is non-degenerate if and only if $\phi(d)$ is non-degenerate.

\begin{definition}
Let ${\cal B}$ be an affine basis of a Delaunay polytope $D$; then

(i) If $rank\,D=1$, then $D$ is called {\em extreme}.

(ii) If all sub-faces of $F(d_{\cal B})$ are degenerate, then $D$ is called {\em maximal}.

\end{definition}
Above definition of maximality is independent of the choice of affine 
basis ${\cal B}$, since, by above Proposition, the linear equivalence
between two faces preserves the non-degeneracy.
Obviously, any extreme Delaunay polytope is maximal.
We present in Corollary \ref{Maximal_Delaunay_polytopes} the list of
all maximal six-dimensional Delaunay polytopes.

Let $\gamma_n=\{0,1\}^n$ be the vertex-set of the Delaunay polytope of the lattice $\Z^n$ and let
\begin{equation*}
h\gamma_n=\{x\in \gamma_n\mbox{~:~}\sum_{i=1}^nx_i\mbox{~is~even}\}
\end{equation*}
be the vertex-set of the Delaunay polytope with center $c=(1/2,\dots, 1/2)$
(this polytope is called {\em half $n$-cube}) of the root lattice 
\begin{equation*}
D_n=\{x\in\Z^n\mbox{~:~}\sum_{i=1}^nx_i\mbox{~is~even}\}\;.
\end{equation*}

\begin{proposition}
The $n$-cube $\gamma_n$ and the half $n$-cube $h\gamma_n$ have rank $n$ and are maximal Delaunay polytopes.
\end{proposition}
\proof One can define the following distance functions on $\gamma_n$ and $h\gamma_n$:
\begin{equation*}
\begin{array}{rcl}
d_i:\gamma_n\times \gamma_n&\rightarrow&\R\\
(x,x')&\mapsto&(x_i-x'_i)^2.
\end{array}
\end{equation*}
If $d$ is the distance function of the $n$-dimensional Delaunay polytope $\gamma_n$ or $h\gamma_n$, then this distance function is expressed as
\begin{equation*}
d=\sum_{i=1}^n \lambda_i d_i\mbox{~with~}\lambda_i>0,
\end{equation*}
which proves that the rank of $\gamma_n$ and of $h\gamma_n$ is $n$. Now, if a face $F$ is included in $F(\gamma_n)$ or $F(h\gamma_n)$, then one of the coefficients $\lambda_i$ becomes zero and the face $F$ corresponds to $\gamma_{n-1}$ or $h\gamma_{n-1}$, which is of lower dimension $n-1$. So, the Delaunay polytopes $\gamma_n$ and $h\gamma_n$ are maximal.

\section{The case of dimension six}\label{Section-DimensionSix}
Call {\em cut cone} and denote by $CUT_{n+1}$ the cone generated by all {\em cuts} $\delta_S\in \R^N$ (where $S$ is a subset of $\{0,\dots ,n\}$), defined by 
\begin{equation*}
(\delta_S)_{ij}=1\mbox{~if~}\vert S\cap \{i,j\}\vert =1\mbox{~and~} (\delta_S)_{ij}=0 \mbox{,~otherwise}.
\end{equation*}
Clearly, $H(b)\delta_S=b(S)(1-b(S))$ with $b(S)=\sum_{a\in S} b_a$; this proves that all hypermetric inequalities are valid on $CUT_{n+1}$. So, $CUT_{n+1}\subset HYP_{n+1}$. Moreover, a cut $\delta_S$ is incident to the face, defined by $H(b)$, if and only if $b(S)=0$ or $1$.

The list of $3773$ facets of $HYP_7$ was found by Baranovskii \cite{Ba99} using the method described in \cite{Ba70}, i.e. he found, by hand, that for all other hypermetric vectors $b$, one can express $H(b)$ as a sum of terms $H(b')$ with $b'$ belonging to his list of $3773$ elements. While this result was announced in \cite{Ba99}, the detailed computations were not published. In \cite{DD-hyp}, another method was proposed: if the Baranovskii's list was not complete, then, in our computation \cite{DD-hyp} of the extreme rays of $HYP_7$, we should find some extreme rays, which are not hypermetric. But this was not the case; so, the list is complete.

The list of representents of $14$ orbits of facets is given below:
\begin{center}
\begin{tabular}{|p{6cm}|p{6cm}|}
\hline
\hline
$b^1=(1,1,-1,0,0,0,0)$&  $b^2=(1,1,1,-1,-1,0,0)$\\\hline
\multicolumn{2}{|l|}{$b^3=(1,1,1,1,-1,-2,0)$, \mbox{~~}$b^4=(2,1,1,-1,-1,-1,0)$}\\\hline
$b^5=(1,1,1,1,-1,-1,-1)$& $b^6=(2,1,1,1,-1,-1,-2)$\\\hline
\multicolumn{2}{|l|}{$b^7=(2,2,1,-1,-1,-1,-1)$, \mbox{~~}$b^8=(1,1,1,1,1,-2,-2)$}\\\hline
\multicolumn{2}{|l|}{$b^9=(3,1,1,-1,-1,-1,-1)$, \mbox{~~}$b^{10}=(1,1,1,1,1,-1,-3)$}\\\hline
\multicolumn{2}{|l|}{$b^{11}=(2,2,1,1,-1,-1,-3)$, \mbox{~~}$b^{12}=(3,1,1,1,-1,-2,-2)$}\\\hline
\multicolumn{2}{|l|}{$b^{13}=(3,2,1,-1,-1,-1,-2)$, \mbox{~~}$b^{14}=(2,1,1,1,1,-2,-3)$}\\
\hline
\hline
\end{tabular}
\end{center}

It gives the total of $3773$ inequalities. The first ten orbits are the orbits of hypermetric facets of the cut cone $CUT_7$; first four of them come as {\em $0$-extension} of facets of the cone $HYP_6=CUT_6$ (see \cite{DL}, Chapter $7$). Last four orbits consist of some $19$-dimensional simplex-faces of $CUT_7$, becoming $20$-dimensional, i.e. simplex-facets in $HYP_7$.

Using Remark \ref{The-Simple-Case-Of-Geometrical-Equivalence}, we obtain that the list of $14$ orbits of facets fall into nine equivalence classes $b^1$, $b^2$, $\{b^3,b^4\}$, $b^5$, $b^6$, $\{b^7,b^8\}$, $\{b^9, b^{10}\}$, $\{b^{11},b^{12}\}$, $\{b^{13}, b^{14}\}$. So, there are nine combinatorial types of six-dimensional Delaunay polytopes of rank $20$ (i.e. repartitioning polytopes).

In \cite{DL} p.~229 another notion, called {\it switching by root} of $b$, is
defined: if $A\subset \{0, \dots, n\}$ and $b(A)=0$, then 
define $b^A$ by $b^A_i=-b_i$ if $i\in A$ and $b^A_i=b_i$ if $i\notin A$.
It is proved that the switching by root of a facet-defining vector of 
$HYP_{n+1}$ is again a facet-defining vector.
Three following vectors define facets of $HYP_7$, which are switching 
by root equivalent:
\begin{equation*}
b^1=(2,2,1,-1,-1,-1,-1), b^2=(2,-2,1,1,1,-1,-1), b^3=(-2,-2,1,1,1,1,1).
\end{equation*}
On the other hand, $b^1$ and $b^3$ are geometrically equivalent by Remark \ref{The-Simple-Case-Of-Geometrical-Equivalence}, while $b^1$ and $b^2$ are not geometrically equivalent.

Remind that $E_6$, $E_7$, $E_8$ are {\em root } lattices defined by
\begin{equation*}
\begin{array}{c}
E_6=\{x\in E_7\mbox{~:~}x_1+x_2=0\},\,\,E_7=\{x\in E_8\mbox{~:~}x_1+\dots+x_8=0\},\\
E_8=\{x\in\R^8\mbox{~:~}x\in \Z^8\cup (\frac{1}{2}+\Z)^8\mbox{~and~}\sum_{i}x_i\in 2\Z\}\;\;.
\end{array}
\end{equation*}

Unique type of Delaunay polytope of $E_6$ is called 
{\it Schlafli polytope} (see \cite{Cox63}) and denoted by $Sch$. 
Its skeleton graph is $27$-vertex (strongly regular) graph, called 
the {\em Schlafli graph}, whose symmetry group has size $51840$
and is isomorphic to the group of isometry, preserving 
the Schlafli polytope. This group is denoted by $Aut(Sch)$.

In \cite{DGL92} were found $26$ orbits of non-cut extreme rays of $HYP_7$ 
by classifying the affine bases of the Schlafli polytope of the root
lattice $E_6$. For every non-cut extreme ray $(\R_+ v)$ of $HYP_7$, there exist
a facet-inducing inequality $f(x)\geq 0$ of $CUT_7$, which is non-hypermetric, 
so that $f(v)<0$. This property establish a bijection between the
$26$ orbits of non-hypermetric facets of the cut cone $CUT_7$ (see \cite{Gcut})
and the $26$ orbits of non-cut extreme rays of $HYP_7$ and proves that 
$HYP_7$ has $29$ orbits of extreme rays:  three orbits of non-zero 
cuts and $26$ orbits coming from $Sch$ (see \cite{DD-hyp}).

\begin{proposition}\label{remarkable-Schlafli}
Let ${\cal B}$ be an affine basis of Schlafli polytope; then

(i) The distance vector $d_{\cal B}$ is incident to $20$ hypermetric faces of $HYP_7$, which are all facets of $HYP_7$.

(ii) If $F$ is a face of $HYP_7$, containing the vector $d_{\cal B}$, then it is non-degenerate and $|Ann(F)|=7+corank(F)$.

\end{proposition}
\proof The Schlafli polytope is six-dimensional and has $27$ vertices. So, for every affine basis ${\cal B}$, the vector $d_{\cal B}$ satisfies $H(b)d_{\cal B}=0$ for $27$ different $b \in Ann(d_{\cal B})$;
we write $Ann(d_{\cal B})-\{e_0,\dots, e_6\}=\{b^1, \dots, b^{20}\}$.
On the other hand, it is known (\cite{DGL92} and \cite{DL} p.~239), that the Schlafli polytope has rank $1$. So, the rank of the matrix $[H(b^1), \dots, H(b^{20})]$ must be ${6+1\choose 2}-1=20$. So, the family $[H(b^i)]_{1\leq i\leq 20}$ is linearly independent. If one of $H(b^i)$ is not a facet of $HYP_7$, then it can be expressed in terms of $[H(b^j)]_{j\not= i}$; this contradicts to linear independence and so, (i) holds.

The extreme ray $d_{\cal B}$ is non-degenerate and is a sub-face of $F$; so, $F$ is also non-degenerate. Every hypermetric face, containing $F$, contains $d_{\cal B}$; so, one has $Ann(F)-\{e_0,\dots,e_6\}=\{{b'}^1, \dots, {b'}^k\}\subset \{b^1, \dots, b^{20}\}$.
The linear independence of the family $[H(b^i)]_{1\leq i\leq 20}$ implies that $k=corank(F)$ and so, (ii) holds.

Above Proposition is not true for the $56$-vertex Gosset polytope (\cite{Cox63}): the Gosset polytope 
has $374$ orbits of affine bases. Each extreme ray, corresponding to an affine basis $\{v_0, \dots, v_7\}$ of the Gosset polytope, is incident to $48$ ($=56-8$) hypermetric faces of $HYP_8$. But amongst these $48$ face-defining inequalities, the number of facets varies from $27$ ($={7+1\choose 2}-1$) to $41$. 
See \cite{DL}, p.~230 for general lower bounds (on the number of 
vertices of a Delaunay polytope) as a function of its rank.

\begin{theorem}\label{ThanksSlava}
All six-dimensional Delaunay polytopes are basic.
\end{theorem}
\proof The simplex is a basic polytope, since there are $7$ vertices and they form an affine basis.

Assume that $D$ is a non-simplicial Delaunay polytope of a lattice $L$ generated by the vectors $w_1$, \dots, $w_6$. Denote by $V$ the volume of the simplex formed by the vectors $0$, $w_1$, \dots, $w_6$. Take a family of $7$ independent vertices ${\cal A}=\{v_0, \dots, v_6\}$ in the vertex-set of $D$ and denote the volume of the corresponding simplex by $V'$.

In \cite{BaRy}, it was proved that the {\em relative volume} $k=\frac{V'}{V}$ is $1$, $2$ or $3$. If $k=1$, then ${\cal A}$ is an affine basis and we are done. Assume now that $k>1$, i.e. that ${\cal A}$ is not an affine basis. Then, there exists a vertex $v$ of $D$, which is written uniquely as $v=\sum_{i=0}^{6} b_i v_i$ with $b$ being fractional and $\sum_{i=0}^6 b_i=1$.

The distance vector $d_{\cal A}$ satisfies $H(b)d_{\cal A}=0$ with $b=(b_0, \dots, b_6)$ being a fractional hypermetric vector. But one can express $H(b)$ as $\sum_{l=1}^{N} \lambda_l H(b^l)$ with $\lambda_l>0$
and $b^l$ being a permutation of one of the following vectors (see \cite{BaRy}):

\begin{center}
\begin{tabular}{|l|l|}
\hline\hline
Case $k=2$                        &Case $k=3$\\
\hline
$\frac{1}{2}(-1,-1, 1,1,1,1,0)$   &$\frac{1}{3}(-1,-1,-1,1,1,2,2)$\\
$\frac{1}{2}(-1,-1,-1,1,1,1,2)$   &$\frac{1}{3}(-1,-1,-1,1,1,1,1)$\\
$\frac{1}{2}(-2,-1,-1,1,1,1,3)$   &$\frac{1}{3}(-2,-1, 1,1,1,1,2)$\\
$\frac{1}{2}(-2,-1,1,1,1,1,1)$    &\\
$\frac{1}{2}(-1,-1,-1,-1,1,2,3)$  &\\
$\frac{1}{2}(-3,-1,1,1,1,1,2)$    &\\
$(-1,1,1,0,0, 0, 0)$              &\\
\hline\hline
\end{tabular}
\end{center}

Since $\lambda_l>0$, one has $H(b^l)d=0$, i.e. the vectors $w'_{l}=\sum_{i=0}^6 b^l_{i} v_i$ with $1\leq l\leq N$ are vertices of the Delaunay polytope $D$.

Since $b$ is fractional, at least one of $b^l$ is fractional, say,
$b^{l_0}$. But all fractional hypermetric vectors of above
Table have one coordinate with absolute value equal to $\frac{1}{k}$,
say, $|b^{l_0}_{i}|=\frac{1}{k}$. So, the family 
$\{v_0, \dots, v_{i-1}, v'_{l_0}, v_{i+1}, \dots, v_6\}$ defines a simplex 
of relative volume $k|b^{l_0}_i|=1$, i.e. it is an affine basis.

\begin{theorem}
The $6241$ combinatorial types of Delaunay polytopes are partitioned by rank in the following way:

{\scriptsize
\begin{center}
\begin{tabular}{|c|c|c|}
\hline
rank     &Nr. in $HYP_7$    &Nr. in $CUT_7$\\\hline
21      &1(simplex)     &0\\
20      &9(repart.)     &1\\
19      &30     &2\\
18      &95     &8\\
17      &233    &28\\
16      &500    &95\\
15      &814    &241\\
14      &1092   &434\\
13      &1145   &527\\
12      &984    &481\\
11      &686    &325\\
10      &417    &183\\
9       &218    &83\\
8       &108    &35\\
7       &52     &13\\
6       &21     &3\\
5       &8      &0\\
4       &4      &0\\
3       &2      &0\\
2       &1      &0\\
1       &1(Schlafli)    &0\\
\hline
\end{tabular}
\end{center}
}
\end{theorem}
\proof The proof is purely computational and the method is described in next Section.

\begin{corollary}\label{Maximal_Delaunay_polytopes}
All maximal six-dimensional Delaunay polytopes are: Schlafli polytope, $6$-cube, half $6$-cube and direct product of half $5$-cube with $1$-cube.
\end{corollary}
\proof This result follows directly from the computation of above Theorem.

\section{Computational methods}\label{Enumeration-Techniques}

Our computation of combinatorial types of six-dimensional Delaunay
polytopes used the face-lattice of $HYP_7$; combinatorial types of
Delaunay polytopes of corank $i+1$ were found from combinatorial types of
Delaunay polytopes of corank $i$. We start from the list of combinatorial
types of corank $1$, i.e. the nine repartitioning polytopes.
The plan of our computation was as follows:
\begin{enumerate}
\item[(i)] Take the list of combinatorial types of Delaunay polytopes of corank $i$.
\item[(ii)] For each of them, find all sub-faces, using our knowledge of facets and extreme rays of $HYP_7$; we obtain faces of corank $i+1$.
\item[(iii)] For every face $F$ of corank $i+1$, find extreme rays $(\R_+ f_i)_{1\leq i\leq N}$, contained in it, and define $d=\sum_{i=1}^N f_i$. The distance vector $d$ is in the interior of $F$; so, using Proposition \ref{BasicReductionResult}, one can test if this ray is non-degenerate or not and this tells us if the face is non-degenerate or not.
\item[(iv)] Find the classes of geometrical equivalence amongst the non-degenerate faces and so, the list of combinatorial types of corank $i+1$.
\end{enumerate}

Above procedure finds all combinatorial types of Delaunay polytopes from corank $1$ (i.e. repartitioning polytopes) till corank $20$ (i.e. the Schlafli polytope). But in order to describe completely our method, we need to precise {\em how} we find the classes under geometrical equivalence.

The first algorithm for the problem of geometrical equivalence is the following: given two faces $F$ and $F'$, find all affine bases of $Ann(F)$ until one finds an equivalence $\phi_{F, F'}$ (see Definition \ref{DefinitionGeometricalEquivalence}).
This method works for corank $1$ or $2$ and was used by Kononenko in the five-dimensional case. But the number of affine bases becomes too important to be workable in corank $3$.

So, one needs another, more efficient method. The first idea is to split the set of faces of $HYP_7$ into two classes: faces, which contain (one or more) Schlafli extreme rays, and those, which are generated only by cuts.

\begin{definition}\label{SubGraphOfSchlafli}
Let $F$ be a face of $HYP_7$, which contains a Schlafli extreme ray $d_{\cal B}$ corresponding to an affine basis ${\cal B}=\{v_0,\dots, v_6\}$. Then, every $b\in Ann(F)$ defines a vertex $v=\sum_{i=0}^6 b_i v_i$ of $Sch$. All such vertices are denoted by ${\cal S}(F,d_{\cal B})$.
\end{definition}

So, every combinatorial type of faces, which contains at least one Schlafli extreme ray, can be interpreted as a set of vertices of $Sch$.

\begin{proposition}
Let $F$ and $F'$ be two faces, which contain at least one Schlafli extreme ray; then it holds:

(i) If $F$ and $F'$ are geometrically equivalent, then for every Schlafli 
extreme ray $d_{\cal B}$ in $F$, there exist an affine basis ${\cal B'}$ 
of $Sch$, such that ${\cal S}(F,d_{\cal B})={\cal S}(F',d_{\cal B'})$.

(ii) If $F$ and $F'$ contain a Schlafli extreme ray $d_{\cal B}$ and 
$d_{\cal B'}$, such that ${\cal S}(F, d_{\cal B})$ is identical to
${\cal S}(F', d_{\cal B'})$ up to an element of $Aut(Sch)$, then 
$F$ and $F'$ are geometrically equivalent.

\end{proposition}
\proof For every distance vector $d_{\cal B}\in F$, one gets an identification of $Ann(F)$ with ${\cal S}(F, d_{\cal B})$.
If $F'$ is geometrically equivalent to $F$, then the vectors $e_0,\dots, e_6$ in $Ann(F')$ are identified with vertices $v'_0, \dots, v'_6$ in ${\cal S}(F,d_{\cal B})$. These vertices form an affine basis ${\cal B'}$ and one has ${\cal S}(F', d_{\cal B'})={\cal S}(F, d_{\cal B})$; so, (i) holds.

Since ${\cal S}(F, d_{\cal B})$ is isomorphic to ${\cal S}(F', d_{\cal B'})$ by an element of $Aut(Sch)$, one can, without loss of generality, assume that they are identical. Now, ${\cal B}$ is identified with vertices $v_0, \dots, v_6$ of ${\cal S}(F, d_{\cal B})$ and ${\cal B'}$ is identified with vertices $v'_0, \dots, v'_6$ of ${\cal S}(F, d_{\cal B})$. Then, the expression of $v'_j$ in terms of $v_i$ determine an affine basis of $Ann(F)$, for which the mapping $\phi_{F,F'}$ is well-defined and bijective.

Above Proposition express the geometrical equivalence in terms of the existence of an element of $Aut(Sch)$ mapping a set of vertices into another set of vertices. Those sets of vertices are identified with corresponding sets of vertices of the Schlafli graph (the Schlafli graph and the Schlafli polytope have the same symmetry group $Aut(Sch)$). So, the problem is expressed in graph-theoretic terms and can be solved, using, for example, the nauty program (\cite{MK}). Therefore, one can build the geometrical equivalence classes.

Now, we extend above method to the case of faces generated by cuts. Let $F$ be a face, generated by cuts $\{\delta(S_i)\}_{1\leq i\leq N}$; we first need to find $Ann(F)$, i.e. all vectors $b\in\Z^{n+1}$ with $\sum_{i=0}^n b_i=1$, having $H(b)\delta(S_i)=0$. 

Those equations can be rewritten as $\sum_{x\in S_i} b_x=x_i$ with $x_i=0$ or $1$, i.e. a linear system in $b$. This linear system has rank $n+1$, because of Proposition below; so, one can find the set $Ann(F)$ for every face, generated by cuts.

\begin{proposition}
Let ${\cal F}$ be a face of $HYP_{n+1}$ generated by cuts $(\delta(S_i))_{1\leq i\leq N}$. Then the following properties are equivalent:

(i) the face $F$ is non-degenerate and

(ii) the linear system, formed by the equations $\sum_{x\in S_i} b_x=0$ and $\sum_{i=0}^n b_i=0$, has solution set $\{0\}$.

\end{proposition}
\proof If the face $F$ is degenerate, then there exist a vertex $v$, which
can be expressed in two different forms 
$v=\sum_{i=0}^n b_i v_i=\sum_{i=0}^n b'_i v_i$. So, after denoting 
$\alpha_i=b_i-b'_i$, one gets $v=\sum_{i=0}^n (b_i+k \alpha_i)v_i$
with $k\in\Z$ and $b+k\alpha$ belongs to $Ann(F)$.
Therefore, one gets $(b+k\alpha)(S_i)=0$ or $1$, and 
$\sum_{i=0}^n b_i+k\alpha_i=1$ for all $k\in \Z$.
This is possible only if $\sum_{x\in S_i} \alpha_x=0$ and 
$\sum_{i=0}^n \alpha_i=0$.

Let the solution set be non-zero, i.e. suppose that one can find an integer-valued non-zero solution $\alpha$. This implies that the vectors $e_0+k\alpha$ belong to $Ann(F)$ for every $k\in\Z$. So, $Ann(F)$ is infinite and $F$ is degenerate.

If we write $Ann(F)=\{b^1, \dots, b^{V}\}$, then every cut $\delta(S_i)\in F$ with $1\leq i\leq N$ defines an Euclidean semi-metric on the set $\{e_0, \dots, e_n\}$; this semi-metric can be uniquely extended to $Ann(F)$ by $\delta_i(b,b')=|b(S_i)-b'(S_i)|$. 

So, to every face, generated by cuts, one can associate a set of semi-metrics on $Ann(F)$, which are, in fact, cut semi-metrics.

A combinatorial type of Delaunay polytope encodes all possible embeddings of $Ann(F)$ into the vertex-set $V$ of a Delaunay polytope of a lattice. These embeddings are completely described by the distance vector $d_{V}$ on their vertices. This distance vector is expressed as $\sum_{i=1}^N \lambda_i \delta_i$ with $\lambda_i>0$. 

Therefore, the combinatorial type of a face, generated by cuts, corresponds to the description of all semi-metrics on this set. This information can be expressed in graph-theoretic terms; so, we can test if two faces, generated by cuts, are isomorphic, using the nauty program (\cite{MK}). So, again one can build the geometrical equivalence classes and our method is completely described.


\begin{thebibliography}{99}

\bibitem[As82]{As}
P. Assouad, {\em Sous-espaces de $L^1$ et in\'egalit\'es hyperm\'etriques}, Compte Rendus de l'Acad\'emie des Sciences de Paris, {\bf 294(A)} (1982) 439--442.

\bibitem[Ba70]{Ba70}
E.P. Baranovskii, {\em Simplexes of $L$-subdivisions of euclidean spaces}, Mathematical Notes, {\bf 10} (1971) 827--834.

\bibitem[Ba99]{Ba99}
E.P. Baranovskii, {\em The conditions for a simplex of $6$-dimensional lattice to be $L$-simplex}, (in Russian) Nauchnyie Trudi Ivanovo State University, Mathematica, {\bf 2} (1999) 18--24.

\bibitem[BK00]{Ko5}
E.P. Baranovski\u\i, P.G. Kononenko, {\em On a method for deriving $L$-polyhedra for $n$-dimensional lattices}, Math. Notes, {\bf 68-5,6} (2000) 704--712.

\bibitem[CS99]{CS}
J.H. Conway and N.J.A. Sloane, {\em Sphere Packings, Lattices and Groups (third edition)}, {\em Grundlehren der mathematischen Wissenschaften}, Springer--Verlag {\bf 290} (1999).

\bibitem[Cox63]{Cox63}
H.S.M. Coxeter, {\em Regular polytopes}, $2^{nd}$ ed., The Macmillan Co. New York 1963.

\bibitem[DGL92]{DGL92}
M. Deza, V.P. Grishukhin, and M. Laurent, {\em Extreme hypermetrics and L-polytopes}, in G.Hal\'asz et al. eds, {\em Sets, Graphs and Numbers, Budapest (Hungary), 1991}, {\bf 60} {\em Colloquia Mathematica Societatis J\'anos Bolyai}, (1992) 157--209.

\bibitem[DGL93]{DGL93}
M.Deza, V.P. Grishukhin, and M. Laurent, {\em The hypermetric cone is polyhedral}, Combinatorica, {\bf 13} (1993) 397--411.

\bibitem[DeLa97]{DL}
M.Deza and M.Laurent, {\em Geometry of cuts and metrics}, Springer--Verlag,
 Berlin 1997.

\bibitem[DD01]{DD-hyp}
M.Dutour and M.Deza, {\em The hypermetric cone on seven vertices}, submitted (2001) and http://il.arXiv.org/abs/math.MG/0108177.

\bibitem[Gr90]{Gcut}
V.P. Grishukhin, {\em All facets of the cut cone $C_n$ for $n=7$ are known}, European Journal of Combinatorics, {\bf 11} (1990) 115--117.

\bibitem[Fe85]{Fe85}
E.S.Fedorov, {\em Elements of the theory of figures} (in Russian)
Imp. Akad. Nauk St.Petersburg 1985 (New edition: Akad. Nauk USSR, 1953).

\bibitem[Ko97]{Ko1}
P.G. Kononenko, {\em Obtaining five-dimensional prime L-polytopes by the
method of layering} (in Russian), Nauchnyie Trudi Ivanovo State University, 
Mathematica, {\bf 1} (1997) 47--55.

\bibitem[Ko98]{Ko2}
P.G. Kononenko, {\em Construction of affine types of L-polytopes of
$5$-dimensional lattices}, deposited in VINITI on 25.11.1998, No 3449-B98, 1--35.

\bibitem[Ko99]{Ko3}
P.G. Kononenko, {\em Affine types of $L$-polytopes of five-dimensional lattices}, (in Russian), Theses (kandidatskaia dissertacia) Ivanovo (1999). 

\bibitem[Ko02]{Ko4}
P.G. Kononenko, {\em Affine types of $L$-polyhedra for five-dimensional lattices}, (in Russian) Mat. Zametki, {\bf 71-3} (2002) 412--430.

\bibitem[MK]{MK}
B.D. McKay, {\em The nauty program}, http://cs.anu.edu.au/people/bdm/nauty/

\bibitem[RyBa98]{BaRy}
S.S. Ryshkov and E.P. Baranovskii, {\em Repartitioning complexes in $n$-dimensional lattices (with full description for $n\leq 6$)}, Voronoi impact on modern science, Book 2, Institute of Mathematics, Kyiv (1998) 115--124.

\bibitem[RyEr87]{ErRyI}
S.S. Ryshkov, R.M. Erdahl, {\em The empty sphere. I}, Canad. J. Math., {\bf 39-4} (1987) 794--824.

\bibitem[RyEr88]{ErRyII}
S.S. Ryshkov, R.M. Erdahl, {\em The empty sphere. II}, Canad. J. Math., {\bf 40-5} (1988) 1058--1073. 

\bibitem[Vo08]{Vo}
G.F. Voronoi, {\em Nouvelles applications des param\`etres continus \`a la th\'eorie des formes quadratiques - Deuxi\`eme m\'emoire}, J. f\"ur die reine und angewandte Mathematik, {\bf 134} (1908) 198-287 and {\bf 136} (1909) 67--178.

\end{thebibliography}
\end{document}